\documentclass[11pt]{amsart}
\usepackage{amscd,amssymb,amsmath}
\newtheorem{theorem}{Theorem}[section]

\theoremstyle{definition}

\newtheorem{conjecture}[theorem]{Conjecture}
\theoremstyle{remark}

\numberwithin{equation}{section}
\begin{document}
\title[Reidemeister numbers of Baumslag - Solitar groups]{Reidemeister number  of any  Automorphism of Baumslag - Solitar group is infinite}
\author{Alexander Fel'shtyn}
\address{ Instytut Matematyki, Uniwersytet Szczecinski,
ul. Wielkopolska 15, 70-451 Szczecin, Poland  
and Boise State University, 1910
University Drive, Boise, Idaho, 83725-155, USA }
\email{felshtyn@diamond.boisestate.edu, felshtyn@mpim-bonn.mpg.de}
\author {Daciberg L. Gon\c{c}alves}
\address{Dept. de Matem\'atica - IME - USP, Caixa Postal 66.281 - CEP 05311-970,
S\~ao Paulo - SP, Brasil}
\email{dlgoncal@ime.usp.br}
\thanks{This work was initiated during the visit of the second author to  Siegen University from September 13 to September 20, 2003.
 The visit was partially supported by a grant
of the ``Projeto tem\'atico Topologia
Alg\'ebrica e Geom\'etrica-FAPESP". The second author would like to thank Professor U. Koschorke for making this visit possible and for the hospitality.\\}

\begin{abstract}
Let $\phi:G \to G$ be a group endomorphism where $G$ is a finitely
generated group of exponential growth, and let  $R(\phi)$ denote
 the number of  $\phi$-conjugacy classes.
Fel'shtyn and Hill \cite{fel-hill} conjectured that if $\phi$ is injective,
then $R(\phi)$ is infinite. This conjecture is true for automorphisms of non-elementary Gromov hyperbolic groups, see \cite{ll} and \cite{fel:1}. It was shown in \cite{gw:2} that the conjecture does not hold in general. Nevertheless in this paper, we show that the conjecture   holds for  the Baumslag-Solitar groups $B(m,n)$, where either $|m|$ or $|n|$ is greater than $1$ and $|m|\ne |n|$. We also show that in the cases where 
$|m|=|n|>1$ or $mn=-1$ the conjecture is true for  automorphisms.
In addition,  we derive  few results  about the coincidence Reidemeister number.

\end{abstract}

\date{\today}
\maketitle

\newcommand{\af}{\alpha}
\newcommand{\et}{\eta}
\newcommand{\ga}{\gamma}
\newcommand{\ta}{\tau}
\newcommand{\ph}{\varphi}
\newcommand{\bt}{\beta}
\newcommand{\lb}{\lambda}
\newcommand{\wh}{\widehat}
\newcommand{\wt}{\widetilde}
\newcommand{\sg}{\sigma}
\newcommand{\om}{\omega}
\newcommand{\cH}{\mathcal H}
\newcommand{\cF}{\mathcal F}
\newcommand{\N}{\mathcal N}
\newcommand{\R}{\mathcal R}
\newcommand{\Ga}{\Gamma}
\newcommand{\cc}{\mathcal C}

\newcommand{\bea} {\begin{eqnarray*}}
\newcommand{\beq} {\begin{equation}}
\newcommand{\bey} {\begin{eqnarray}}
\newcommand{\eea} {\end{eqnarray*}}
\newcommand{\eeq} {\end{equation}}
\newcommand{\eey} {\end{eqnarray}}

\newcommand{\ovl}{\overline}
\newcommand{\vv}{\vspace{4mm}}
\newcommand{\lra}{\longrightarrow}

\bibliography{ref}
\bibliographystyle{amsplain}

\section{Introduction }

J. Nielsen introduced the fixed point  classes of a surface homeomorphism  in \cite{ni}.
 Subsequently, K. Reidemeister \cite{reid:re} developed the algebraic foundation of the  Nielsen fixed point theory for any map of any compact polyhedron. As a result of Reidemeister's work,
the twisted conjugacy  classes of a group  homomorphism were  introduced. It turns out that
the fixed point  classes of a map  can easily be
identified  with the conjugacy classes of
lifting of this map to the universal covering of compact polyhedron,  and conjugacy classes of lifting can   be identified
with the twisted conjugacy  classes of the homomorphism induced on  the
fundamental group of the polyhedron.
Let $G$ be a finitely generated group and let $\phi: G\rightarrow G$ be  an endomorphism.
Two elements $\alpha,\alpha^\prime\in G$ are said to be
$\phi-$$conjugate$ if there exists $\gamma \in G$ with
$
\alpha^\prime=\gamma  \alpha   \phi(\gamma)^{-1}.
$
The number of  $\phi$-conjugacy classes( or twisted conjugacy classes )
 is called the $Reidemeister$ $number$ of an  endomorphism $\phi$, denoted by $R(\phi)$. If $\phi$ is the identity map, then the $\phi$-conjugacy classes are the usual conjugacy classes in the group $G$.
Let $X$  be a connected  compact
polyhedron and $f:X\rightarrow X$  be a continuous map.
The Reidemeister number $R(f)$, which  is simply the cardinality of
the set of  $\phi$-conjugacy classes where $\phi=f_{\#}$ is the
induced homomorphism on the fundamental group, is relevant for the study of fixed points of $f$  in the presence of the fundamental group. In fact
the finiteness of Reidemeister number  plays an important r\^ole.
See for  example  \cite{wong2}, \cite{gw:3},
\cite{fel-hill}, \cite{FelTroObzo} and the introduction of \cite{gw:2}.

A current  important  problem  concerns  obtaining
a twisted analogue of the celebrated Burnside-Frobenius theorem
\cite{fel-hill, FelshB, FelTroVer, polyc, FelTroObzo}.
 For this purpose
it is important to describe the class of groups $G$, such that $R(\phi)=\infty$
for any automorphism $\phi:G\to G$. 
A. Felshtyn and R. Hill \cite{fel-hill} made  first attempts to localize this class of groups.

Later  it was proved in \cite{fel:1,ll} that the
non-elementary Gromov hyperbolic groups belong to this class.
Furthermore, using the co-Hofian property, it was  shown  in \cite{fel:1}  that, if in addition $G$ is torsion-free and freely indecomposable,  then $R(\phi)$ is infinite for every injective endomorphism  $\phi$.
This result gives supportive evidence to a conjecture of
\cite{fel-hill} which states that $R(\phi)=\infty$ if $\phi$ is an injective endomorphism   of  a finitely generated torsion-free group $G$  with exponential growth.

This conjecture was shown  to be false in general. In \cite {gw:2} were
constructed  automorphisms $\phi: G\to G$ on certain finitely
generated torsion-free exponential growth groups $G$ that are not Gromov
hyperbolic with $R(\phi)<\infty$.

In the present  paper we study this  problem for a family of finitely generated groups
 which have exponential growth but are not Gromov hyperbolic.  These  are the
Baumslag-Solitar groups, which we now define.
 Being  indexed by pairs of integer numbers different from zero, they  have the following
presentation:   $$B(m,n)=\langle a,b: a^{-1}b^ma=b^n \rangle, m,n\ne 0 .$$
\noindent The present work  extends substantially in several directions the preliminary  results obtained in \cite{fg}, and  simplifies  some of the  proofs. 
 
The  family of the Baumslag-Solitar groups has different features from the one given in \cite{gw:2}, which is a family  of  groups which are metabelian  having as  the kernel the group
${\Bbb Z}^n$. Hence  they contain a subgroup isomorphic to ${\Bbb Z} + {\Bbb Z}$. In the case of Baumslag-Solitar groups this happens if, and only if,  $m=n$.
For $m=n=1$ the group $B(1,1)={\Bbb Z}+{\Bbb Z}$ does not have exponential growth and it is also known that there are automorphisms $\phi: B(1,1) \to B(1,1)$ with  $R(\phi)<\infty$.
For more details about these groups $B(m,n)$ see \cite{b-s, fa-mo1}.

Some results in this work could be obtained by means of the classification of some of the endomorphisms   of a  Baumslag-Solitar group (for those, see \cite{GHMP} and \cite{koch}).
We use only one direct consequence of the main result of  \cite{koch} which concerns  injective homomorphisms.

Our main results are:\\
\noindent{\bf Theorem } For any injective endomorphism of $B(m,n)$ where 
$|n|\ne |m|$ and $|nm|\ne 0$, the Reidemeister number is infinite.
For  any automorphism of $B(m,n)$ where $0<|m|=|n|$ and $mn\ne 1$, the  Reidemeister number  is also  infinite.

This result summarizes  the results of    Theorems  3.4,  4.4, , 5.4, 6.4  and Proposition 5.1 for the various values of $m$ and $n$.

\noindent{\bf Theorem 7.1 }  The coincidence  Reidemeister number    is infinite for any pair of injective endomorphisms  of the group
$B(m,n)$,  where $|n|\ne |m|$ and $|nm|\ne 0$.

We do not know if Theorems  5.4 and 6.4 are  also true for injective homomorphisms. See Remarks 5.5 and 6.5.

We say that  a group $G$ has \emph{property } $R_\infty$ if any of  its automorphisms  $\phi$ has $R(\phi)=\infty$.
After the preprint( in arXiv: math.GR-0405590) of this   article was circulate, 
it was proved that the  following groups   have 
 \emph{property } $R_\infty$ :
(1) generalized Baumslag-Solitar groups, that is, finitely generated groups
which act on a tree with all edge and vertex stabilizers infinite cyclic
\cite{LevittBaums}; \ 
(2) lamplighter groups $Z_n \wr Z$ iff $2|n$ or $3|n$ \cite{gowon1}; \ (3)
the solvable generalization $\Gamma$ of $BS(1,n)$ given by the short exact sequence
$1 \rightarrow Z[\frac{1}{n}] \rightarrow \Gamma \rightarrow Z^k \rightarrow 1$
as well as any group quasi-isometric to $\Gamma$ \cite{TabWong}; \ 
(4) groups which are quasi-isometric to the Generalized Baumslag-Solitar  groups  \cite{TabWong2} (while this property is not a quasi-isometry invariant); \  (5)
saturated weakly branch groups( including the Grigorchuk group and the Gupta-Sidki group) \cite{FelLeoTro}; \  (6)  the 
R.~Thompson's   groups \cite{ThomsonRinf}; \ (7) some finitely generated nilpotent groups
of arbitrary Hirsch length \cite{gowon2}.

 We would like to complete the introduction with the following
conjecture.

\begin{conjecture}
Any relatively  hyperbolic   group  has \emph{property } $R_\infty$.
In particular,  any Kleinian group has \emph{property } $R_\infty$.
\end{conjecture}

This  paper is organized into six sections besides this one. In section 2, we make some simple reduction of the problem to certain cases and  develop some preliminaries about the Reidemeister classes of a pair of homomorphisms  between short exact sequences. In section 3, we  study the case $B(\pm 1,n)$ for $|n|>1$, with  main result  Theorem 3.4.
 In section 4, we consider the cases $B(m,n)$ for $1<|m|\ne |n|>1$, with main result   Theorem 4.4. In section 5, we consider the cases $B(m,-m)$ for $|m|>0$, with  main results Propositon 5.1 and   Theorem 5.4. In section 6 we consider the  cases $B(m,m)$ for $|m|>1$, with main result  Theorem 6.4. In section 7 we derive  few results about the coincidence Reidemeister number, with main result Theorem 7.1 \\

{\bf Acknowledgments} The authors would like to thank G. Levitt for his helpful comments improving an  earlier  version of this manuscript.
The first  author would like to thank
 B. Bowditch, T. Januszkiewicz,  M. Kapovich  and E. Troitsky  for stimulating discussions and comments. The second author would like to express his  thanks to D.  Kochloukova for very helpful  discussions.\\
The first  author also would  like  to thank the Max-Planck-Institute f\"ur Mathematik, Bonn for its kind hospitality and support.

This article is dedicated to the memory of Sasha Reznikov.

\vskip 1.0cm

\section{Generalities and Preliminaries}

\vskip 1.0cm

In this section we first  describe few elementary properties of the groups 
$B(m,n)$ in order to reduce our problem to certain cases. Then  we recall some facts about the  Reidemeister classes of a pair of homomorphisms  of a short exact sequence. Recall that  {\it a group $G$ has \emph{property } $R_\infty$ if any of  its automorphisms  $\phi$ has $R(\phi)=\infty$.}

Recall that the Baumslag-Solitar groups are indexed by pairs of integer numbers different from zero  and  they  have the following
presentation:    
  $$B(m,n)=\langle a,b: a^{-1}b^ma=b^n \rangle, m,n\ne 0 .$$ 

 The first observation is that for $m=n=1$ this group is ${\Bbb Z} + {\Bbb Z}$. It is well known that this group does not have exponential growth and there are automorphisms $\phi:{\Bbb Z} + {\Bbb Z} \to {\Bbb Z} + {\Bbb Z}$ with finite Reidemeister number. So  ${\Bbb Z} + {\Bbb Z}$ does not have property $R_\infty$. 

The second observation is that $B(m,n)$ is isomorphic to $B(-m,-n)$. It suffices to see that the relations  $a^{-1}b^ma=b^n$ and 
 $a^{-1}b^{-m}a=b^{-n}$ each one generates  the same normal subgroup, since one relation is the inverse of the other.

The last observation is that  $B(m,n)$ is isomorphic to $B(n,m)$. Suppose that  $B(m,n)=\langle a,b: a^{-1}b^ma=b^n \rangle, m,n\ne 0 $ and 
$B(n,m)=\langle c,d: c^{-1}d^nc=d^m \rangle, m,n\ne 0 .$ The map which sends $a \to c^{-1}$ and $b \to d$ extends to an isomorphism of the two groups.

Based on the above, we will consider only the groups $B(r,s), \ rs\ne 0,1$ and we can show

\noindent{\bf Proposition 2.0} Each group $B(r,s), \ rs\ne 0,1 $ is isomorphic to some $B(m,n),$ where  $m,n$ satisfy $0<m\leq |n|$ and $n\ne 1$.  

 So we will divide the problem into 4 cases. Case 1) is when $1=m<|n|$; Case 2) is when $1<m<|n|$; Case 3) when $0<m=-n$; Case 4) when $1<m=n$ 

\bigskip 


The set of the Reidemeister classes of a pair of homomorphisms will be denoted by $R[\ ,\ ]$ and  the number of such classes by
$R(\ , \ )$. When  the two sequences are the same and one of the homomorphisms  is the identity, then we have the usual Reidemeister classes and Reidemeister number. The main reference for this section is   \cite{go:nil1} and more details can be found there.

Let us consider a diagram of two short exact sequences of groups and
maps between these two sequences:
\beq
\begin{array}{ccccccccc}
1 & \to & H_1 & \stackrel{i_1}{\to} & G_1 & \stackrel{p_1}{\to} & Q_1 &
\to & 1 \\
&&&&&&&& \\
& & f'\downarrow\downarrow g' & & f\downarrow\downarrow g & & \ovl f\downarrow
\downarrow \ovl g & & \\
&&&&&&&& \\
1 & \to & H_2 & \stackrel{i_2}{\to} & G_2 & \stackrel{p_2}{\to} & Q_2
& \to & 1
\end{array}
\eeq
where $f' = f|_{H_1}$, $g' = g|_{H_1}$.

We recall that the set of the  Reidemeister classes $R[f_1, f_2]$ relative to
homomorphisms $f_1, f_2 : K \to \pi$ is the set of  the equivalence classes of
elements of $\pi$ given by the following relation: $\alpha \sim f_2(\ta)
\alpha f_1(\ta)^{-1}$ for $\alpha \in \pi$ and $\tau \in K$.

 The diagram $(2.1)$   provides  maps between sets \[ R[f',g'] \stackrel{\wh i_2}{\to}R[f,g] \stackrel{\wh p_2}{\to} R[\ovl f, \ovl g]\] where the last map is clearly surjective. An  obvious consequence of this  fact will be used to solve some of the cases that we will  discuss, and that will appear below as Corollary 2.2. For the remaining cases we need further   information about the above sequence and we will use   Corollary 2.4.

Proposition 1.2 in    \cite{go:nil1} says

\noindent {\bf Proposition 2.1}  Given the diagram $(2.1)$ we have a short sequence of sets
\[ R[f',g'] \stackrel{\wh i_2}{\to}R[f,g] \stackrel{\wh p_2}{\to}
R[\ovl f, \ovl g]\]
where $\wh p_2$ is surjective and $\wh p_2^{-1}[1] = {\rm im}\,(\wh i_2),$
where $1$ is the identity element of $Q_2$.

An immediate consequence  is

\noindent {\bf Corollary 2.2} If $ R(\ovl f, \ovl g)$
is infinite,  then $ R(f,g)$ is also infinite.

\begin{proof} Since  $\wh p_2$ is surjective the result follows.
\end{proof}

In order to study the injectivity of the map $\wh i_2$,  for each element
$\ovl\alpha \in Q_2$ let $H_2(\ovl\alpha) = p_2^{-1}(\ovl\alpha)$,  $C_{\ovl\alpha} =
\{\ovl\ta \in Q_1|\ovl g(\ovl\ta)\ovl\alpha\ovl f(\ovl\ta^{-1}) = \ovl\alpha\}$ 
and let $R_{\ovl\alpha}[f',g']$ be the set of equivalence classes of elements
of $H_2(\ovl\alpha)$ given by the equivalence relation $\beta \sim g(\ta)\beta f
(\ta^{-1}),$ where $\beta \in H_2(\ovl\alpha)$ and $\ta \in p_1^{-1}(C_{\ovl\alpha})$.
Finally, let $R[f_{\ovl\af}, g_{\ovl\alpha}]$ be the set of equivalence classes
of elements of $H_2(\ovl\alpha)$ given by the relation $\beta \sim g(\ta)\beta f
(\ta^{-1})$, where $\beta \in H_2(\ovl\alpha)$ and $\ta \in G_1$.

Proposition 1.2 in    \cite{go:nil1} says

\noindent{\bf Proposition 2.3}
Two classes of $R(f_{\ovl\alpha}, g_{\ovl\alpha})$  represent the same
class of $ R(f,g)$ if and only if  they belong to the same orbit by
the action of $C_{\ovl\alpha}$. Further the isotropy subgroup of this
action at an element $[\beta]$ is $G_{[\beta]} = p_1(C_{\beta}) \subset
C_{\ovl\alpha}$ where $\beta \in [\beta]$.

An immediate consequence of this result is

\noindent {\bf Corollary 2.4} If  $C_{\ovl\alpha}$ is finite and $R(f_{\ovl\alpha}, g_{\ovl\alpha})$ is infinite for some $\alpha$,  then
 $ R(f,g)$ is also infinite. In particular,  this is the case if $Q_2$ is finite.

\begin{proof} The orbits of the action of $C_{\ovl\alpha}$ on  $R[f_{\ovl\alpha}, g_{\ovl\alpha}]$ are finite. So we have an infinite number of orbits. The last part is an easy consequence of  the first part.
\end{proof}

 \vskip 1.0cm

\section{ THE CASES $B(m,n)$, $1=|m|<|n|$}

\vskip 1.0cm

From section 2 the cases in this section  reduce to  Case 1), namely  
$B(1,n)$ for $1<|n|$.  Let $|n|>1$ and $B(1,n)=\langle a,b: a^{-1}ba=b^n, n>1 \rangle $.\\
Recall from \cite{fa-mo1} that the Baumslag-Solitar groups $B(1,n)$
are finitely generated solvable groups which are not
virtually nilpotent. These  groups have
exponential growth \cite{harpe},   and  they are not Gromov hyperbolic.
Furthermore, those groups are metabelian and  torsion  free.

Consider the homomorphisms $|\ \ |_a: B(1,n) \longrightarrow {\Bbb Z} $ which associates for each word $w\in B(1,n)$ the sum of the exponents of $a$ in the word. It is easy to see that this is a well defined map into $\Bbb Z$ which is surjective.

\noindent {\bf Proposition 3.1} We have a short exact sequence
$$0 \longrightarrow K\longrightarrow B(1,n) \stackrel{|\ \ |_a}\longrightarrow   {\Bbb Z} \longrightarrow 1,$$
\noindent  where $K$,  the kernel of the
map  $ |\ \ |_a$, is the set of  the elements which have  the sum of the powers  of $a$ equal to zero. Furthermore, $B(1,n)=K \rtimes {\Bbb Z}$ (semi-direct product).

\begin{proof} The first part is clear. The second part follows because ${\Bbb Z}$ is free, so the sequence splits.
\end{proof}

\noindent {\bf Proposition 3.2} The kernel $K$ coincides with  $ N\langle b\rangle, $  the   normalizer of $\langle b \rangle$ in $B(1,n)$.\\

\begin{proof}
We have $ N\langle b\rangle\subset K $. But the quotient $B/N\langle b\rangle
$ has the following presentation: $\bar a^{-1}\bar b \bar a= \bar b^n, \bar b=1$. Therefore this group is isomorphic to ${\Bbb Z}$ and the natural projection coincides with the map $|\ \ |_a$ under the obvious identification of ${\Bbb Z}$ with  $B/N\langle b\rangle$. Consider the commutative diagram \\

$
\begin{matrix}
0 {\rightarrow}    N\langle b\rangle {\rightarrow} & B(1,n) {\rightarrow}  &   B/N\langle b\rangle
{\rightarrow} & 1  \cr
 \downarrow   & \downarrow & \downarrow  \cr
  0 {\rightarrow} K {\rightarrow}  & B(1,n) {\rightarrow} & {\Bbb Z} {\rightarrow} & 1  \end{matrix}
$

\noindent where the last vertical map is an isomorphism. From the well-known
 five Lemma the result follows.
\end{proof}

 The groups $B(1,n)$ are metabelian. Let $\epsilon$ be the sign of $n$.  We recall the result that  $B(1,n)$
is isomorphic to ${\Bbb Z}[1/|n|]\rtimes_{\theta} {\Bbb Z}$ where the action of ${\Bbb Z}$ on ${\Bbb Z}[1/|n|]$ is given by 
$\theta(1)(x)=x/n^{\epsilon}$. To see this, first observe that  the map defined by   $\phi(a)=(0,1)$ and  $\phi(b)=(1,0)$ extends to a  unique homomorphism $\phi:B \to  {\Bbb Z}[1/|n|]\rtimes {\Bbb Z}$ which is clearly surjective. It suffices to show that this homomorphism is injective. Consider  a word $w=a^{r_1}b^{s_1}...a^{r_t}b^{s_t}$ such that $r_1+...+r_t=0$. Thus   $w\in K$, and using   the relation of the group this word is equivalent to $b^{s_1/n^{\epsilon r_1}}b^{s_2/n^{\epsilon(r_1+r_2)}}....b^{s_{t-1}/n^{\epsilon(r_1+...+r_{t-1})}}.b^{s_t}$.
If we apply $\phi$ to this element, which belongs to the kernel of $\phi$, we obtain that the sum of the powers  ${s_1/n^{r_1}}+{s_2/n^{\epsilon(r_1+r_2)}}+....+{s_{t-1}/n^{\epsilon(r_1+...+r_{t-1})}}+{s_t}$ is zero. But  this means that $w$ is the trivial element, hence $\phi$ restricted to $K$ is injective. Therefore the
result follows.

\noindent{\bf Proposition 3.3} Any homomorphism $\phi: B(1,n) \to B(1,n)$ is a homomorphism of the short exact sequence given in Proposition 3.2.

\begin{proof} Let $\bar \phi$ be the  homomorphism induced by $\phi$ on the abelianization  of $B(1,n)$. The abelianization  of  $B(1,n)$, denoted by   $B(1,n)_{ab}$, is isomorphic to $\Bbb Z_{|n-1|} +\Bbb Z$. The torsion elements of $B(1,n)_{ab}$ form a subgroup isomorphic to $\Bbb Z_{|n-1|}$ which is invariant under any homomorphism. The preimage of this subgroup under the projection to the abelianization  $B(1,n) \to B(1,n)_{ab}$ is  exactly the subgroup $N(b)$, i.e., the elements represented by words where the sum of the powers of $a$ is zero. So it follows that $N(b)$ is mapped into $N(b)$.
\end{proof}
\noindent{\bf Theorem 3.4} For any injective homomorphism of $B(1,n)$  the  Reidemeister number  is infinite.

\begin{proof} Let $\phi$ be an injective endomorphism. By Proposition 3.3 it is an endomorphism of the short exact sequence given by Proposition 3.2.  The induced homomorphism  on the quotient   is a non-trivial  endomorphism of $\Bbb Z$.  Otherwise we would have an injective homomorphism from the non-abelian group $B(1,n)$ into the abelian group $K$.  If the induced endomorphism
 $\bar \phi$ is the identity,
by Corollary     2.2   the number of Reidemeister classes is  infinite and the result follows. So,  assume   that $ \bar \phi$ is
multiplication by $k\ne 0,1$ and we will get a contradiction. Now we claim that there is no injective 
endomorphism of $B(1,n)$ such that the induced homomorphism  on the quotient is multiplication
by $k$ with $k\ne 0,1$. When we apply the homomorphism
$\phi$ to the relation $a^{-1}ba=b^n$, using the isomorphism $B(1,n) \to \Bbb Z[1/n]\rtimes \Bbb Z$
  we obtain: $a^{-k}\phi(b)a^k= (n^k\phi(b),0)=(n\phi(b),0)$, which implies that either  
$n^{1-k}=1$ or $\phi(b)=0$. Since $\phi(b)\ne 0$ and  $n$ is neither 1 or -1 we get a contradiction and  the result follows.
\end{proof}

\noindent{\bf Remark 3.5} From the proof above we conclude that any injective homomorphism $\varphi:B(1,n)\to B(1,n)$ has the property that it induces the identity on the quotient ${\Bbb Z}$ given by the short exact sequence in Propositon 3.2. This fact will be used to study coincidence Reidemeister classes in section 7.

\vskip 1.0cm

\section { THE CASE $B(m,n)$, $1<|m|\ne|n|>1$}

\vskip 1.0cm

From section 2 the cases in this section  reduce to  Case 2), namely  $(m,n)$ for $1<m<|n|$. The groups in this section are more complicated than the ones in the previous
section. Nevertheless in order to obtain the results we will use a  similar procedure to  the one in the previous section.
Let $1<m<|n| $ and $B(m,n)=\langle a,b: a^{-1}b^ma=b^n \rangle $.
Recall that such groups are non-virtually solvable.

Consider the homomorphism $|\ \ |_a: B(m,n) \longrightarrow {\Bbb Z} $ which associates to each word $w\in B(m,n)$ the sum of the powers of $a$ in the word. It is easy to see that this is a well defined homomorphism into ${\Bbb Z}$ which is surjective.

\noindent {\bf Proposition 4.1} We have a short exact sequence $$0 \longrightarrow K\longrightarrow B(m,n) \longrightarrow {\Bbb Z} \longrightarrow 1,$$ where $K$, the kernel of the
map  $ |\ \ |_a$, is the set of  the elements which have  the sum of the powers  of $a$ equals to zero.  Furthermore, $B(m,n)=K \rtimes {\Bbb Z}$ is a semi-direct product where the action is given with respect to some  fixed section.

\begin{proof} The first part is clear. The second part follows because ${\Bbb Z}$ is free, so the sequence splits. Since the kernel $K$ is not abelian,  the action is defined with respect to a specific section (see \cite{br}).
\end{proof}

\noindent{\bf Proposition 4.2} The kernel $K$ coincides with  $ N\langle b\rangle $ which is the   normalizer of $\langle b \rangle$ in $B(m,n)$.\\

\begin{proof} Similar to Proposition 3.2.
\end{proof}

\noindent{\bf Proposition 4.3} Any homomorphism $\phi: B(m,n) \to B(m,n)$ is a homomorphism of the short exact sequence  given in Proposition 4.1.

\begin{proof} Let $\bar \phi$ be the homomorphism induced by $\phi$ on the abelianized of $B(m,n)$. The abelianized of  $B(m,n)$, denoted by   $B(m,n)_{ab}$, is isomorphic to ${\Bbb Z}_{|n-m|} +{\Bbb Z}$. The torsion elements of $B(m,n)_{ab}$ form a subgroup isomorphic to ${\Bbb Z}_{|n-m|}$ which is invariant under any homomorphism. The preimage of this subgroup under the projection to the abelianized $B(m,n) \to B(m,n)_{ab}$ is  exactly the subgroup $N(b)$, i.e., the elements represented by words where the sum of the powers of $a$ is zero. So it follows that $N(b)$ is mapped into $N(b)$.
\end{proof}

In order to have a homomorphism $\phi$ of $B(m,n)$ which  has  finite Reidemeister number, the induced map on the quotient ${\Bbb Z}$ must be
different from the identity by the same argument used in the proof of Theorem 3.4. 

Now we will give a presentation of the group $K$.  The  group $K$ is generated by the elements
$g_i= a^{-i} b a^i  \ i\in {\Bbb Z}$ which  satisfy the following relations:
$\{1= a^{-j} (a^{-1} b^m a b^{-n}) a^j = g_{j+1}^m g_j^{-n}\}$   for all integers $j$.
This presentation is a consequence of the Bass-Serre theory, see \cite{co},
Theorem 27,  page 211.

Now we will prove the main result of this section. Denote by  $K_{ab}$ the abelianization  of $K$.

\noindent{\bf Theorem 4.4}  For any injective homomorphism of $B(m,n)$  the  Reidemeister number  is infinite.

\begin{proof}   Let us  consider the short exact sequence,
obtained from the short exact sequence given in  Proposition 4.1,  by taking the quotient
 with the commutators subgroup of $K$, i.e.
$$0 \longrightarrow K_{ab}\longrightarrow B(m,n)/[K,K] \longrightarrow {\Bbb Z} \longrightarrow 1.$$
 Thus  we obtain a short exact sequence where the kernel $K_{ab}$ is abelian. From the
presentation of $K$ we obtain a presentation of $K_{ab}$ given as follows:
 It is generated by the elements
$g_i,  \ i\in Z$,  which  satisfy the following relations:
$\{1= g_{j+1}^m g_j^{-n}, g_ig_j=g_jg_i \}$   for all integers $i,j$. This presentation  is the same
as the quotient of the free abelian group generated by the elements 
$g_i, \ i\in {\Bbb Z}$ (so the direct sum of ${\Bbb Z}'s$ indexed by ${\Bbb Z}$), modulo the subgroup generated by the relations 
$\{1= g_{j+1}^m g_j^{-n} \}$. Thus  an element can be regarded as an equivalence class of
a sequence of integers indexed by ${\Bbb Z}$, where the elements of the sequence are 
zero but  a finite number. By abuse of notation
we also  denote  by $\phi$  the induced endomorphism on $ B(m,n)/[K,K]$. 

Let  $\phi(a)=a^{k}\theta$ for $\theta\in K_{ab}$ and $k\ne 1$. Recall that if $k=1$
it follows immediately that the Reidemeister number is infinite.  Since the kernel of the extension is abelian, after applying  $\phi$ to the relation 
$a^{-1}b^ma=b^n$  we obtain \\
 $$\theta^{-1}a^{-k}\phi((b)^m)a^{k}\theta=a^{-k}\phi(b)^ma^{k}= \phi(b^n)=\phi(b)^n.$$
\noindent From the main result of \cite{koch}, the element $\phi(b)$ is a conjugate of a power of $b$, i.e.,  $\phi(b)=\gamma b^r\gamma^{-1}$ for some $r\ne 0$. In the abelianization  the element $\gamma b^r\gamma^{-1}$ is the same as the element $a^sba^{-s}$ for some integer $s$.  So any power of 
 $\phi(b)$ with exponent different from zero is a nontrivial element. Now we take both sides of the equation above to the power $m^k$. If $k=0$  it  follows immediately that $m=n$. Let us take $k>1$. After applying  the relation several times we obtain 
$$ a^{-k}\phi(b)^{mm^k}a^{k}= \phi(b)^{mn^{k}}=\phi(b)^{nm^{k}}.$$
Therefore it  follows that $mn^{k}=nm{k}$ or $ n^{k-1}=m^{k-1}$. If $n$ is positive,  since $k\ne 1$, then the   only  solutions are $m=n$,  which is a contradiction. If $n$ is negative then the only  solutions are,
$m=n$ or $m=-n$ and $k$ even. In either case we get a contradiction.  The case where $k<0$ is similar and the result follows.
 \end{proof}

\noindent{\bf Remark 4.5} The  Proposition 4.1 certainly holds for $n=-m$. If we apply the proof of the Theorem 4.4 above
for the group $B(m,-m)$  we can conclude that any injective homomorphism 
$\varphi:B(m,-m)\to B(m,-m)$ has the property that it  induces the multiplication by an odd number on the quotient ${\Bbb Z}$, where  
${\Bbb Z}$ is given by the short exact sequence used in the proof of the Theorem 4.4.

\vskip 1.0cm

\section {THE CASE $B(m,-m)$, $1\leq |m|$}

\vskip 1.0cm

From section 2 the cases in this section  reduce to  Case 3), namely 
$B(m,-m)$ for $0<m$. 

We will start with the group $B(1,-1)$. In this case it   has been proved in (see \cite{gowon2}) that this group has the $R_{\infty}$ property.  For sake of completeness we write another proof (which is also known to the second author of \cite{gowon2}) where the techniques  is more close to the ones 
 used in this work. The group $B(1,-1)$ is isomorphic to the fundamental group of the Klein bottle.

\noindent {\bf Proposition 5.1} For any automorphism $\phi$ of  ${\Bbb Z}\rtimes {\Bbb Z}$ the Reidemeister number is   infinite.

\begin{proof}  Let us consider the short exact sequence 
$$0\to {\Bbb Z}\to {\Bbb Z}\rtimes {\Bbb Z} \to {\Bbb Z},$$
\noindent where the inclusion ${\Bbb Z}\to {\Bbb Z}\rtimes {\Bbb Z}$ sends  $1  \to x$. It is well know that ${\Bbb Z}$ is characteristic in ${\Bbb Z}\rtimes {\Bbb Z}$, so any homomorphism  
$\varphi:  {\Bbb Z}\rtimes {\Bbb Z}\to {\Bbb Z}\rtimes {\Bbb Z}$ induces a homomorphism of short exact sequence. Let $\varphi $ be an automorphism. Then the induced  automorhism on the quotient $\bar \varphi:{\Bbb Z} \to {\Bbb Z}$ is either the identity or minus the identity. In the first case we have that the Reidemeister number of $\varphi$ is infinite and the result follows. So let us assume that $\bar \varphi$ is $-id$. The induced map on the fiber $\varphi'$ is also either the identity or minus the identity. In either case, in order to compute the Reidemeister number of $\varphi$, by means of the formula given in \cite{gw:2}, Lemma 2.1 we need to consider the homomorphism given by the composition  of $\varphi'$ with the conjugation by $y$, which is the multiplication by -1, i.e. $-\varphi'$. So either $\varphi'$ or $-\varphi'$ is  the identity. Again  by the formula given in \cite{gw:2}, Lemma 2.1,  the result follows. 

\end{proof}

 The above result is not true for injective homomorphisms. Take for example the homomorphism defined by $\varphi(x)=x^2, \varphi(y)=y^3$. It is an injective homomorphism and $R(\varphi)$ is 4. 

 From now on let $1<m$.  The groups $B(m,-m)$, in contrast with   others Baumslag-Solitar groups already considered,   have  subgroups isomorphic to  ${\Bbb Z}\rtimes {\Bbb Z}=B(1,-1)$,  the fundamental group of the Klein bottle.

It is straightforward  to verify that   Propositions 4.1, 4.2 and 4.3  are also true for $m=-n$ (this is not the case for Proposition 4.3 when $m=n$).  So we have a short exact sequence 
  $$0 \longrightarrow K\longrightarrow B(m,-m) \longrightarrow {\Bbb Z} \longrightarrow 1,$$  where $K$ is  the kernel of the
map  $ |\ \ |_a$. This kernel coincides with the normal subgroup generated by the element $b$ and any homomorphism $\phi: B(m,-m) \to B(m,-m)$ is a homomorphism of the short exact sequence. Denote by $\bar \phi$ the induced homomorphism  on the quotient ${\Bbb Z}$ and by  $\phi':K \to K$ the restriction of $\phi$.  Our proof have some similarities
with the proof for the group $B(1,-1)$. 

\noindent{\bf Proposition 5.2}  Given any automorphism  $\phi:B(m,-m)\to B(m,-m)$, then the induced automorphism  on the quotient $\bar \phi$ is either the identity or minus the identity. In the former case we have 
$R(\phi)$ infinite.

\begin{proof} Follows immediatly from Corollary 2.2.

\end{proof}

\noindent{\bf Proposition 5.3}  Given any automorphism  $\phi:B(m,-m)\to B(m,-m)$ such that the  induced homomorphism on the quotient 
$\bar \phi:{\Bbb Z} \to {\Bbb Z}$ is multiplication by $-1$, either the automorphism $\phi'$ or the automorphism $\tau_a\circ \phi'$, where  
$\tau_a$ is the conjugation by the element $a$, have infinite Reidemeister number.

\begin{proof} From \cite{koch} $\phi'(b)=a^ib^{\epsilon}a^{-i}$ for some integer $i$,  where $\epsilon$ is either $1$ or $-1$ since we have an automorphism. For the  other automorphism  we have  
$\tau_a\circ \phi'(b)=a^{i+1}b^{\epsilon}a^{-i-1}$. We certainly have
either $\epsilon=(-1)^i$ or   $\epsilon=(-1)^{i+1}$. Let $\varphi$ be
$\phi'$ or  $\tau_a\circ \phi'$ according to  $\epsilon=(-1)^i$ or   
$\epsilon=(-1)^{i+1}$, respectively. A presentation of the group $K$ was given in section 4 before Theorem 4.4. Consider the extra relation in $K$ 
given by  $b=ab^{-1}a^{-1}$. Then it follows that  the quotient group is ${\Bbb Z}$ and   the automorphism $\varphi$ induces a homomorphism on the quotient which agrees with the identity. So it  follows that the Reidemeister number of $\varphi$ is infinite.
\end{proof}



\noindent{\bf Theorem 5.4} For any automorphism of $B(m,-m)$  the  Reidemeister number  is infinite.

\begin{proof} Let $\phi: B(m,-m)\to B(m,-m)$ an automorphism. From Proposition 5.2 we can assume that $\bar \phi$ is multiplication by $-1$. From Propositon 5.3 we know that either  $\phi'$ or  $\tau_a\circ \phi'$
has infinite  Reidemeister number. By the formula given in \cite{gw:2}, Lemma 2.1,  the result follows. 
\end{proof}

\noindent{\bf Remark 5.5} We do not know an example of an injective homomorphism on $B(m,-m),$ for $m>1$, which has finite  Reidemeister number.


\vskip 1.0cm

\section {THE CASE $B(m,m)$, $|m| > 1$}

\vskip 1.0cm

From section 2 the cases in this section  reduce to  Case 4), namely 
$B(m,m)$ for $1<m$. The proof of this case is not similar to the previous cases.

As we noted  before,   if $m=1$,  the group is ${\Bbb Z} + {\Bbb Z}$. Then there are   automorphisms  of the group  which  have a finite number of Reidemeister classes. For $m>1$, in order to study its  automorphisms, we describe the groups $B(m,m)$ as certain extensions. Finally  we show that any automorphism of this group has  infinite Reidemeister number.

These groups, in contrast with  the  Baumslag-Solitar groups already considered,   have  subgroups isomorphic to  ${\Bbb Z} + {\Bbb Z}$. We remark that for $n=2$ this is not the fundamental group of the Klein bottle. There is a surjection from $B(2,2)$ onto the fundamental group of the Klein bottle.

We start by describing these groups.
Let  $|\ \ |_b: B(m,m) \to {\Bbb Z} $ be the homomorphism which associates to a word the sum of the powers of $b$ which appears in the word. This is a 
well defined surjective homomorphism and we have

\noindent{\bf Proposition 6.1}  There is a  splitting  short exact sequence:
$$0  \to F \to B(m,m) \to {\Bbb Z} \to 1,$$
\noindent where  $F$ is the free group on $m$ generators $x_1,...,x_m$ and
the last map is $|\ \ |_b$.
Further, the action of the generator $1\in Z$ is the automorphism of $F$ which,  for $j<m$, sends $x_j$ to $x_{j+1}$ and $x_m$ to $x_1$.

\begin{proof} Let $F\rtimes {\Bbb Z}$ be the semi-direct product of $F$ by ${\Bbb Z},$ where $F$ is the free group on  $x_1,...,x_m$ and the action is given by the automorphism of $F$ which, for $j<m$,  maps  $x_j$ to $x_{j+1}$ and $x_m$ to $x_1$.
We will show that $B(m,m)$ is isomorphic to  $F\rtimes {\Bbb Z}$. For this  consider the map $\psi: \{a,b\} \to F\rtimes {\Bbb Z}$ which sends $a$ to $x_1$ and $b$ to $1\in Z$. This map extends to a homomorphism $B(m,m) \to F\rtimes Z$, which we also denote by $\psi$,
 since the relation which defines the group $B(m,m)$ is preserved by the map. Also $\psi$  is a homomorphism  of short exact sequences.
The map restricted to the kernel of $| \ \  |_b$ is surjective to the free group  $F$. Also the kernel admits a set of generators with cardinality $n$. So the map restricted  to the kernel is an isomorphism and the result follows.

 \end{proof}

  Proposition 6.1  above shows that the groups $B(m,m)$ are policyclic.
For more about the Reidemeister number of these groups see \cite {fgw}.

\noindent{\bf Proposition 6.2}  The center of $B(m,m)$ is the subgroup generated by $b^m$. Moreover,   any injective homomorphism $\phi:B(m,m) \to B(m,m)$ leaves the center invariant.

\begin{proof} For the first part, from Proposition 6.1, we know that $B(m,m)$ is of the form  $F\rtimes {\Bbb Z}$. Let $(w,b^r)\in F\rtimes {\Bbb Z}$ be in the center and $(v,1)\in F\rtimes {\Bbb Z}$ where $v$ is an arbitrary element of $F$. We have $(w,b^r)(v,1)=(w\cdot b^r(v),b^r)$ and $(v,1)(w,b^r)=(v\cdot w,b^r)$. We can assume that $w$ is a word written in the reduced form which starts with $x_i^{m_i}$, for some $1 \leq i \leq m$. Let $r_0$ be the integer,  $0 \leq r_0 \leq m-1$,  congruent to $r$ mod $m$. Now we consider three  cases:\\
 
(1)  $r_0=0$. Take $v=x_{i+1}$ if $i<m$ and $v=x_1$ if $i=m$. We claim that $w.b^r(v)\ne v.w$, so the elements do not commute.
 To see that they do not commute  observe first that $v.w$ is in the reduced form. If $w.b^r(v)$ is not reduced they cannot be equal. If it is reduced, also they can not be equal either,  since they start with different letters. The  argument above does not work if $w=1$, but this is the case where the element is in the center.\\

(2)   $r_0\ne 0$ and $w\ne 1$. Take $v=x_i^{m_i}$. Again  $v.w$ is in the reduced form  starting with $x_i^{2m_i}$. If $w\cdot b^r(v)$ is not reduced they cannot be equal. If it is reduced, also they cannot be equal either,  since they  start with different powers of $x_i$, even if the word contains only one letter,  since $b^r(v)$ is not a power of $x_i$ ($r$ is not congruent to $0$ mod $m$). \\

(3)  $r_0\ne 0$ and $w=1$. Then $r=km+r_0$,  and from the relation in the group it  follows that  $a^{-1}b^ra=a^{-1}b^{km+r_0}ra=b^{km}a^{-1}b^{r_0}a$. But $a^{-1}b^{r_0}a=b^{r_0}$ implies $b^{r_0}ab^{-r_0}=a$,  which in terms of the notation in Proposition 5.1 means $x_1=x_{r_0}$, which is a contradiction. So the result follows.

For the second part we have to show that $\phi(b^m)$ is in the center.
Since $\phi$ is injective, from the main result of \cite{koch}, the element 
$\phi(b)$ is  conjugated to  a power of $b$, i.e.,  $\phi(b)=\gamma b^r\gamma^{-1}$ for some $r\ne 0$. Therefore   $\phi(b^m)=\gamma (b^r)^m\gamma^{-1}=bmr$ and the result follows. 

\end{proof}

Next  we consider the group which is the quotient of $F\rtimes {\Bbb Z}$ by the center, where the center   is the subgroup $<b^m>$. This quotient is isomorphic to
$F\rtimes {\Bbb Z}_n$ and we denote the image of the  generator $b$ in ${\Bbb Z}$  by $\bar b$ in ${\Bbb Z}_m$. 

\noindent{\bf Proposition 6.3}  Any  automorphism of the group $F\rtimes {\Bbb Z}_m$   has   infinite  Reidemeister number. 

\begin{proof} We know that $F$ is the free group on the letters $x_1,..., x_m$ and let 
$\theta :F\rtimes {\Bbb Z} \to {\Bbb Z}_n$ be the homomorphism  defined by $\theta(x_i)=1$ and $\theta(\bar b)=0$. The kernel of this homomorphism defines a subgroup   of $F\rtimes {\Bbb Z}$ of index $m$ which is isomorphic to   $F'\rtimes {\Bbb Z}_n$,  where $F'$ is the kernel of the homomorphism $\theta$ restricted to $F$.
 Now we claim that $F'$ is invariant with respect to any  homomorphism, i.e,  $F'$ is characteristic. Let $(w,\bar 1)$ be an arbitrary element of the subgroup $F'$ with $w\ne 1$. First observe that $\theta(\phi(x_i))=\theta(\phi(x_1))$, for all $i$. This follows by induction since  $x_{i+1}=\bar b.x_i. \bar b^{-1}$,   
$\theta(\phi(x_{i+1}))=\theta(\phi(\bar b)).\theta(\phi(x_i)).\theta( \phi(\bar b^{-1}))=\theta(\phi(x_i)).$ Therefore
$\theta(\phi(w,\bar 1))=\theta((w,\bar 1))\theta(\phi(x_1))$ and  the subgroup $F'$ is invariant. Therefore the   automorphism $\phi$ provides an automorphism  of the short exact sequence $$0 \to F' \to F\rtimes {\Bbb Z}_m \to {\Bbb Z}_m+{\Bbb Z}_m \to 0$$
where the restriction to the kernel is an automorphism of a free group of finite rank. Hence,  by the Corollary 2.4 the result follows. 
\end{proof}

Now we proof the main result.

\noindent{\bf Theorem 6.4} Any automorphism   $\phi$ of $B(m,m)$ has
an infinite Reidemeister number.

\begin{proof} Any  automorphism $\phi$, from Proposition 5.2,  induces 
an automorphism on $F\rtimes {\Bbb Z}_m$, which we  denote by $\bar \phi$. In order to prove that $\phi$ has an infinite  Reidemeister number,  it suffices to show the same statement  for  $\bar \phi$.
By  Proposition 5.3 the statement  is true for $\bar \phi$,  so the Theorem  follows.
\end{proof}

{\bf Remark 6.5}  Proposition 6.3 and Theorem  6.4 use only  Proposition 6.2  for automorphisms, which, under this assumption, its  second part of the Proposition 6.2  becomes obvious. Nevertheless, using Proposition 6.2 as stated,  it is not difficult to see that  Proposition 6.3 and  Theorem 6.4  can be extended  for injective homomorphisms  if one knows the result for injective homomorphisms of a free group of finite rank. However  this  is still an open question.

\section {Coincidence Reidemeister classes}

For a pair of homomorphisms $\phi,\psi:G \to G$ one can ask   similarly    when a  pair of homomorphisms $(\phi,\psi)$   has infinite  coincidence Reidemeister number. If one of the homomorphisms, let us say $\phi$, is an automorphism,  then  the problem is equivalent to the classical problem for the homomorphism 
$\phi^{-1}\circ \psi$. So we can apply all the  results above.  There are many cases which can be obtained from the case of one homomorphism.  Theorems 3.4,  4.4, 6.2,  and 6.4   can be generalized to   coincidence as follows:

\noindent{\bf Theorem 7.1 }  The coincidence  Reidemeister number    is infinite for any pair of injective endomorphisms  of the group
$B(m,n)$,  where $|n|\ne |m|$ and $|nm|\ne 0$  .\\

 \begin{proof} For the cases in question, we have proved that an injective  homomorphism induces a  homomorphism of the short exact sequences given by Propositions 3.1 and   4.1,  depending on the values of $m$ and $n$, respectively.  
Further in any of the cases above, by the proof of the Theorems 3.4,  4.4,  6.2 and 6.4, we have that the induced  homomorphisms $\bar \phi$ and $\bar \psi$ on  the quotients are the identity   on $\Bbb Z$. So the pair  $(\bar \phi,\bar \psi)$ has  infinite  coincidence
 Reidemeister number  and the result follows from Corollary 2.2.
\end{proof}

The extension of Theorem 7.1 for the groups $B(m,m)$,   will  follow  if the same result  is true   for a pair of injective   homomophisms of a free group of finite rank. But, as pointed out in Remark 6.5,  this is not known even if one of the homomorphisms  is the identity.


\end{document}